\documentclass[graybox]{svmult}
\usepackage{amssymb,amsmath, amsfonts, amsbsy, amsbsy, latexsym, color, epsfig, graphicx}
\usepackage{mathptmx}       % selects Times Roman as basic font
\usepackage{helvet}         % selects Helvetica as sans-serif font
\usepackage{courier}        % selects Courier as typewriter font
\usepackage{type1cm}        % activate if the above 3 fonts are
\usepackage{makeidx}         % allows index generation
\usepackage{graphicx}        % standard LaTeX graphics tool
\usepackage{multicol}        % used for the two-column index
\usepackage[bottom]{footmisc}% places footnotes at page bottom

\makeindex             % used for the subject index
\newcommand{\bitem}{\begin{itemize}}
\newcommand{\eitem}{\end{itemize}}
\newcommand{\benum}{\begin{enumerate}}
\newcommand{\eenum}{\end{enumerate}}
\newcommand{\beq}{\begin{equation}}
\newcommand{\eeq}{\end{equation}}
\newcommand{\ip}[2]{\langle#1,#2\rangle}

\newcommand{\norm}[1]{\|#1\|}

\newcommand{\Sp}{\mbox{supp}}
\newcommand{\R}{\mathbb{R}}

\newcommand{\Z}{\mathbb{Z}}

\def\eps{\varepsilon}

\def\cC{{\mathcal{C}}}

\def\cE{{\mathcal{E}}}
\def\cQ{{{Q}}}
\def\cR{{\mathcal{R}}}

\def\cS{{{S}}}

\def\cSH{{{S}{H}}}
\def\ZZ{\mathbb{Z}}
\def\RR{\mathbb{R}}
\def\ZZ{\mathbb{Z}}

\newcommand{\bZ}{{\mathbb Z}}
\newcommand{\bR}{{\mathbb R}}

\def\cC{{\mathcal{C}}}

\def\cR{{\mathcal{R}}}

\newcommand{\wql}[1]{{\color{black}{#1}}}
\newcommand{\wq}[1]{{\color{black}{#1}}}
\newcommand{\gk}[1]{{\color{black}{#1}}}

\newcommand{\gggk}[1]{{\color{black}{#1}}}
\newcommand{\jk}[1]{{\color{black}{#1}}}

\begin{document}

\title*{Shearlets on Bounded Domains}
\author{Gitta Kutyniok and Wang-Q Lim}
\institute{Gitta Kutyniok \at Institute of Mathematics, University of Osnabr\"uck, 49069 Osnabr\"uck, Germany, \email{kutyniok@math.uos.de}
\and Wang-Q Lim \at Institute of Mathematics, University of Osnabr\"uck, 49069 Osnabr\"uck, Germany, \email{wlim@math.uos.de}}
\maketitle

\abstract{Shearlet systems have so far been only considered as a means to analyze $L^2$-functions defined on $\RR^2$, which
exhibit curvilinear singularities. However, in applications such as image processing or numerical solvers of partial
differential equations the function to be
analyzed or efficiently encoded is typically defined on a non-rectangular shaped bounded domain.
Motivated by these applications, in this paper, we first introduce a novel model for cartoon-like images defined on a bounded domain. We then prove that
compactly supported shearlet frames satisfying some weak decay and smoothness conditions, when
orthogonally projected onto the bounded domain, do provide (almost)
optimally sparse approximations of elements belonging to this model class.
}

%*******************************************************************************************
%*******************************************************************************************
%*******************************************************************************************

\section{Introduction}
\label{sec:1}

It is by now well accepted that $L^2$-functions \gggk{supported} on the unit cube which are $C^2$ except for a $C^2$ discontinuity curve are
a suitable model for images which are governed by edges. Of all directional representation systems which
provide optimally sparse approximations of this model class, shearlet systems have distinguished themselves by the
fact that they are the only system which provides a unified treatment of the continuum and digital setting, thereby making
them particularly useful for both theoretical considerations as well as applications.
However, most applications concern sparse approximations of functions on bounded domains, for instance, a
numerical solver of a transport dominated equation could seek a solution on a polygonal shaped area. This calls for shearlet
systems which are adapted to bounded domains while still providing optimally sparse expansions.

In this paper, we therefore consider the following questions:
\begin{enumerate}
\item[(I)] Which is a suitable model for a function on a bounded domain with curvilinear singularities?
\item[(II)] What is the `correct' definition of a shearlet system for a bounded domain?
\item[(III)] Do these shearlet systems provide optimally sparse approximations of the model functions introduced in (I)?
\end{enumerate}
In the sequel we will indeed provide a complete answer to those questions. These results push the door open for the
usability of shearlet systems in all areas where 2D functions on bounded domains require efficient
encoding.

%*******************************************************************************************
%*******************************************************************************************

\subsection{Optimally Sparse Approximations of Cartoon-Like Images}
\label{subsec:cartoonmodel}

The first complete model of \jk{cartoon-like images} has been introduced in \cite{CD04}, the basic idea being that a closed $C^2$ curve
separates smooth -- in the sense of $C^2$ -- functions. For the precise definition, we let $\rho : [0,2\pi] \to \RR^+$
be a \gk{$C^2$} radius function and define the set $B$ by
\begin{equation}\label{eq:curve}
B = \{x \in \RR^2 : \|x\|_2 \le \rho(\theta), \: x = (\|x\|_2,\theta) \mbox{ in polar coordinates}\},
\end{equation}
where
\begin{equation}\label{eq:curvebound}
\sup|\rho^{''}(\theta)| \leq \nu, \quad \rho \leq \rho_0 < 1.
\end{equation}
This allows \jk{us} to introduce $STAR^2(\nu)$, a class of sets $B$ with $C^2$ boundaries $\partial B$ and
curvature bounded by $\nu$, as well as $\cE^2(\nu)$, a class of cartoon-like images.
%\gggk{In the sequel,} we will denote
%the space of twice continuously differentiable functions on \gggk{$\RR^2$ with compact support in $[0,1]^2$} by $C^2([0,1]^2)$.

\begin{definition}[\cite{CD04}] \label{def:def1}
For $\nu > 0$, the set $STAR^2(\nu)$ is defined to be the set of all $B \subset [0,1]^2$ such that $B$
is a translate of a set obeying \eqref{eq:curve} and \eqref{eq:curvebound}. Further, $\cE^2(\nu)$
denotes the set of functions $f$ \gggk{on $\RR^2$ with compact support in $[0,1]^2$} of the form
\[
f = f_0 + f_1 \chi_{B},
\]
where \gggk{$B \in STAR^2(\nu)$ and $f_0,f_1 \in C^2(\RR^2)$ with compact support in $[0,1]^2$ as well as
$\sum_{|\alpha| \leq 2}\norm{D^{\alpha}f_i}_{\infty} \leq 1$ for each $i=0,1$.}
\end{definition}

In \cite{Don01}, Donoho  proved that the optimal rate which can be achieved under some restrictions on the representation
system as well as on the selection procedure of the approximating coefficients is
\[
\norm{f-f_N}_2^2 \le C \cdot N^{-2} \quad \mbox{as } N \to \infty,
\]
where $f_N$ \gk{is the best $N$-term approximation}.

%*******************************************************************************************
%*******************************************************************************************

\subsection{Shortcomings of this Cartoon-Like Model Class}

The first shortcoming of this model is the assumption that the discontinuity curve is $C^2$. Think, for instance, of an image, which
pictures a building. Then the frames of the windows separate the dark interior of the windows from the presumably light color of the wall,
however this frame is far from being $C^2$. Hence, a much more natural assumption would be to assume that the discontinuity curve is
piecewise $C^2$.

The second shortcoming consists in the fact that the function is implicitly assumed to vanish on the boundary of $[0,1]^2$. More precisely,
even if the function $f = f_0 + f_1 \chi_{B}$ is non-zero on a section of positive measure of the boundary $\partial B$, this situation
is not particularly treated at all. However, reminding ourselves of the very careful boundary treatment in the theory of partial
differential equations, this
situation should be paid close attention. Thus, a very natural approach to a careful handling of the boundary in a model for cartoon-like images
seems to consist in regarding the boundary as a singularity curve itself.

The third and last shortcoming is the shape of the \gggk{support} $[0,1]^2$ of this model. Typically, in real-world situations the domain of
2D data can be very different from being a rectangle, and even a polygonal-shape model might not necessarily be sufficient. Examples
to support this claim can be found, for instance, in fluid dynamics, where the flow can be supported on \gk{variously} shaped domains.
In this regard, a suitable model situation seems to be to allow the boundary to consist of any piecewise $C^2$ curve.

%*******************************************************************************************
%*******************************************************************************************

\subsection{Our Model for Cartoon-Like Images on Bounded Domains}
\label{subsec:ourmodel}

\begin{figure}[h]
\begin{center}
\includegraphics[height=1.25in]{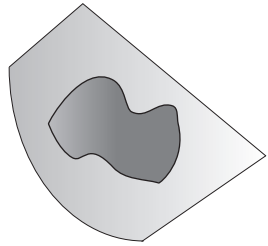}
\end{center}
\caption{Example of a function $f$ belonging to our model class  $\cE^2_{\nu,L}(\Omega)$ of cartoon-like images on bounded domains.}
\label{fig:Cartoon}
\end{figure}

The model for cartoon-like images on bounded domains, which we now define, will indeed take all considerations from the previous subsection
into account. For an illustration, we refer to Figure \ref{fig:Cartoon}.

We first introduce $STAR^2(\nu,L)$, a class of sets $B$ with now piecewise $C^2$ boundaries $\partial B$ and
curvature on each piece bounded by $\nu$. This will serve us for both modeling the bounded domain as well as
modeling the discontinuity curve.  For this, let $L \in \Z^+$ denote the number of $C^2$ pieces and let $\nu > 0$
be an upper estimate for the curvature on each piece.
Then $B \in STAR^2(\nu,L)$, if $B$ is a bounded subset of $[0,1]^2$ \wq{whose boundary
$\partial B$ is a simple closed curve} and
%\ggk{Here I suggested something else; see the file from June 24! Why did you prefer
%this version?}\\ \wql{I think that $\sup |\rho''| \leq \nu$ needs to be changed to something like $$\max_{i=1,\dots,L}\sup |\rho_i''| \leq \nu$$ in your definition.
%Notice that $\rho$ is not $C^2$ but piecewise $C^2$ smooth. On the other hand, the number of smooth sub-curves $\rho_i$, $L$ does not appear in your definition.
%Finally, I thought that it might be better to use a parametrization in terms of Cartesian coordinates (in stead of polar) to indicate that
%our shearlet setting is more convenient for the Cartesian coordinates rather than polar - also this parametrization directly appears in the proof, which makes the proof more clear.
%What do you think ?  }
\[
\partial B = \bigcup_{i=1}^{L} \rho_i
\]
where each curve $\rho_i$ is parameterized by either $x_1 = E_i(x_2)$ or $x_2 = E_i(x_1)$ and $E_i \in C^2([a_i,b_i])$
such that
\[
\max_{i=1,\dots,L}\jk{\max_{[a_i,b_i]}} |E^{''}_i| \leq \nu.
\]

This allows us to introduce a model class of cartoon-like images in bounded domains. \gggk{In accordance with
modeling functions on bounded domains, we now consider functions defined on $[0,1]^2$; its `true' domain is brought into play
by requiring these functions to be supported on $\Omega \subseteq (0,1)^2$, which we model as piecewise $C^2$ bounded. This
ensures that we treat $\partial \Omega$ as a singularity curve, which would not have been possible when defining
the model on $\Omega$ itself.}

\begin{definition}
For $\nu>0$ and $L \in \Z^+$, let $\Omega, B \in STAR^2(\nu,L)$ be such that $B \subset \Omega^\circ$, where $\Omega^\circ$ denotes
the interior of the set $\Omega$, and $\Omega \subset (0,1)^2$.  Then, $\cE^2_{\nu,L}(\Omega)$ denotes the set of functions
$f$ \gggk{on $[0,1]^2$ with compact support in $\Omega$} of the form
\[
f = f_0 + f_1 \chi_{B},
\]
where \gggk{$f_0,f_1 \in C^2([0,1]^2)$ with compact support in $\Omega$} and $\sum_{|\alpha| \leq 2}\norm{D^{\alpha}f_i}_{\infty} \leq 1$ for each $i=0,1$.
\end{definition}

%\gggk{We remark that in the sequel we might also sometimes regard a function $f \in \cE^2_{\nu,L}(\Omega)$ as a function
%defined on $\Omega$ itself.}

Later it will become important to analyze the points on boundaries of sets in $STAR^2(\nu,L)$, in which the boundary is
not $C^2$. For these points, we will employ the following customarily used notion.

\begin{definition}
For $\nu>0$ and $L \in \Z^+$, let $B \in STAR^2(\nu,L)$. Then a point $x_0 \in \partial B$ will be called a {\em corner point},
if $\partial B$ is not $C^2$ in $x_0$.
\end{definition}

Since the model $\cE^2_{\nu,L}(\Omega)$, while containing the previous model $\cE^2_{\nu}$ as a special case, is considerably more complicated, we
would like to make the reader aware of the fact that it is now not clear at all whether the optimal approximation rate is
still
\[
\norm{f-f_N}_2^2 \le C \cdot N^{-2} \quad \mbox{as } N \to \infty.
\]

%*******************************************************************************************
%*******************************************************************************************

\subsection{Review of Shearlets}

The directional representation system of {\em shearlets} has recently emerged -- a first introduction
dates back to 2005 in \cite{LLKW05} -- and rapidly gained
attention due to the fact that, in contrast to other proposed directional representation systems,
shearlets provide a unified treatment of the continuum and digital world similar to wavelets.
We refer to, e.g., \cite{GKL06,KL09} for the continuum theory, \cite{DKS08,ELL08,Lim09} for the digital
theory, and \cite{GLL09,DK10} for recent applications. Shearlets are scaled according to a parabolic
scaling law and exhibit directionality by parameterizing slope by shearing, the later being the secret
which allows the aforementioned unified treatment in contrast to rotation. Thus shearlets are associated
with three parameters: scale, orientation, and position. A precise definition will be given in
Section \ref{sec:2}.

A few months ago, the theory of shearlets focussed entirely on band-limited generators although
precise spatial localization is evidently highly desirable for, e.g., edge detection. Recently,
motivated by this desideratum, compactly supported shearlets were studied by Kittipoom and the
two authors. It was shown that a large class of compactly supported shearlets generate a frame
for $L^2(\RR^2)$ with controllable frame bounds alongside with several explicit constructions
\cite{KKL10a}. By the two authors it was then proven in \cite{KL10} that a large class of these
compactly supported shearlet frames does in fact provide (almost) optimally sparse approximations of
functions in $\cE^2_{\nu}$ in the sense of
\[
\norm{f-f_N}_2^2 \le C \cdot N^{-2} \cdot (\log N)^{3} \quad \mbox{as } N \to \infty.
\]
It should be mentioned that although the optimal rate is not completely achieved, the $\log$-factor is
typically considered negligible compared to the $N^{-2}$-factor, wherefore the term `almost optimal'
has been adopted into the language.

%*******************************************************************************************
%*******************************************************************************************

\subsection{Surprising Result}
\label{subsec:surprise}

We now aim to analyze the ability of shearlets to sparsely approximate elements of the previously introduced model for
cartoon-like images on bounded domains, $\cE^2_{\nu,L}(\Omega)$. For this, we first need to
define shearlet systems for functions in $L^2(\Omega)$. Assume we are given a (compactly supported) shearlet
frame for $L^2(\RR^2)$. The most crude approach to transform this into a shearlet system defined on $L^2(\Omega)$,
where $\Omega \in STAR^2(\nu,L)$,
is to just truncate each element at the boundary of $\Omega$. Since it is well known in classical frame theory
that the orthogonal projection of a frame onto a subspace does not change the frame bounds (cf. \cite{Chr03}),
this procedure will result in a (compactly supported) shearlet frame for $L^2(\Omega)$ with the same
frame bounds as before.

We now apply this procedure to the family of compactly supported shearlet frames for $L^2(\RR^2)$, which yielded
(almost) optimally sparse approximations of functions in $\cE^2_{\nu}$ (see \cite[Thm. 1.3]{KL10}). The main
result of this paper then proves that the resulting family of shearlet frames \gggk{ -- now regarded as a system
on $[0,1]^2$ with compact support in $\Omega$ -- } again provides (almost)
optimally sparse approximations now of elements from our model of cartoon-like images on bounded domains
$\cE^2_{\nu,L}(\Omega)$ in the sense of
\[
\norm{f-f_N}_2^2 \le C \cdot N^{-2} \cdot (\log N)^{3} \quad \mbox{as } N \to \infty.
\]
The precise statement is phrased in Theorem \ref{theo:main} in Section \ref{sec:3}.

This result is quite surprising in two ways:
\bitem
\item {\em Surprise 1}. Regarding a $\log$-factor as negligible -- a customarily taken viewpoint --, the previous result shows that
even for our much more sophisticated model of cartoon-like images on bounded domains the {\em same} optimal
sparse approximation rate as for the simple model detailed in Subsection \ref{subsec:cartoonmodel} can be achieved.
This is even more surprising taking into account that our model contains point singularities at the corner
points of the singularity curves. Naively, one would expect that these should worsen the
approximation rate. However, observing that `not too many' shearlets intersect these `sparsely occurring' points
unravels this mystery.
\item {\em Surprise 2}. Orthogonally projecting a shearlet system onto the considered boun\-ded domain, thereby
merely truncating it, seems an exceptionally crude approach to derive shearlets for a bounded domain. However,
these `modified' shearlet systems are indeed sufficient to achieve the optimal rate and no sophisticated
adaptions are required, which is of significance for deriving fast algorithmic realizations.
\eitem

%*******************************************************************************************
%*******************************************************************************************

\subsection{Main Contributions}

\gggk{The main contributions of this paper are two-fold. Firstly, we introduce $\cE^2_{\nu,L}(\Omega)$ as a suitable model
for a function on a bounded domain with curvilinear singularities. Secondly, we show that the `crude' approach towards a
shearlet system on a bounded domain by simply orthogonally projecting still provides optimally sparse approximations
of elements belonging to our model class $\cE^2_{\nu,L}(\Omega)$.

We should mention that although not formally stated the idea of one piecewise $C^2$ discontinuity curve in a model for functions on
$\RR^2$ as an extension of Definition \ref{def:def1} is already lurking in \cite{CD04}. Also a brief sketch of proof of
(almost) optimally sparse approximations of curvelets is contained therein. These ideas are however very different from
ours in two aspects.
First of all, our goal is a suitable model for functions on bounded domains exhibiting discontinuity curves and also treating the
boundary of the domain as a singularity curve. And secondly, in this paper we consider compactly supported shearlets
-- hence elements with superior spatial localization properties in contrast to the (band-limited) curvelets -- which
allows an elegant proof of the sparse approximation result in addition to a simplified treatment of the corner points.}

%*******************************************************************************************
%*******************************************************************************************

\subsection{Outline}

In Section \ref{sec:2}, after recalling the definition of shearlet systems, we introduce
shearlet systems on bounded domains, thereby focussing in particular on compactly supported shearlet frames.
The precise statement of our main result is presented in Section \ref{sec:3} together with a road map
to its proof. The proof itself is then carried out in Section \ref{sec:4}. Finally, in Section \ref{sec:5},
we discuss our results and possible extensions of it.

%*******************************************************************************************
%*******************************************************************************************
%*******************************************************************************************

\section{Compactly Supported Shearlets}
\label{sec:2}

We first review the main notions and definitions related to shearlet theory, focussing in
particular on compactly supported generators. \gk{For more details we would like to refer the
interested reader to the survey paper \cite{KLL10}.} Then we present our definition of shearlet systems on a
bounded domain $\Omega \in STAR^2(\nu,L)$.

%*******************************************************************************************
%*******************************************************************************************

\subsection{Compactly Supported Shearlet Frames for $L^2(\RR^2)$}

Shearlets are scaled according to a parabolic scaling law encoded in the {\em parabolic scaling matrices}
$A_{2^j}$ or $\tilde{A}_{2^j}$, $j \in \ZZ$, and exhibit directionality by parameterizing slope encoded
in the {\em shear matrices} $S_k$, $k \in \ZZ$, defined by
\[
A_{2^j} =
\begin{pmatrix}
  2^j & 0\\ 0 & 2^{j/2}
\end{pmatrix}
\qquad\mbox{or}\qquad
\tilde{A}_{2^j} =
\begin{pmatrix}
  2^{j/2} & 0\\ 0 & 2^j
\end{pmatrix}
\]
and
\[
S_k = \begin{pmatrix}
  1 & \wq{k}\\ 0 & 1
\end{pmatrix},
\]
respectively.

We next partition the frequency plane into four cones $\cC_1$ -- $\cC_4$. This allow the introduction of
shearlet systems which treat different slopes equally in contrast to the shearlet group-based approach.
We though wish to mention that historically the shearlet group-based approach was developed first due to
very favorable theoretical properties and it still often serves as a system for developing novel analysis
strategies (see, for instance, \cite{KKL10b}).

\begin{figure}[ht]
\begin{center}
\includegraphics[height=1.4in]{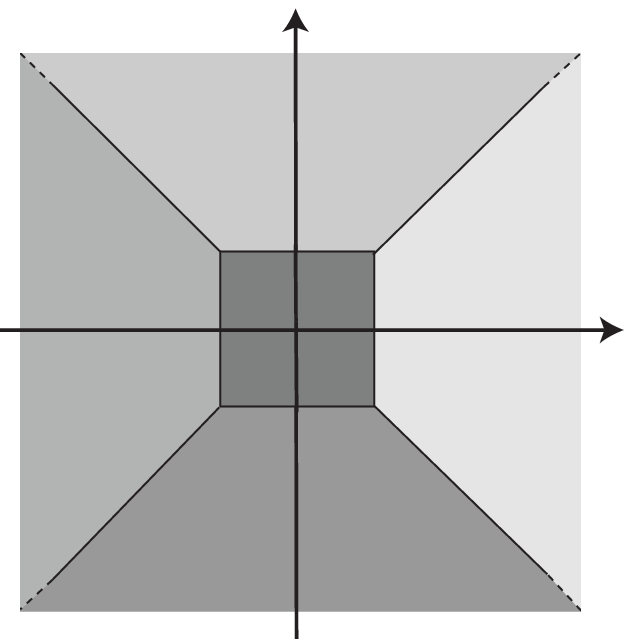}
\put(-33,58){\footnotesize{$\cC_1$}}
\put(-70,80){\footnotesize{$\cC_2$}}
\put(-88,30){\footnotesize{$\cC_3$}}
\put(-50,52){\footnotesize{$\cR$}}
\put(-45,15){\footnotesize{$\cC_4$}}
\put(-60,-17){(a)}
\hspace*{2cm}
\includegraphics[height=1.4in]{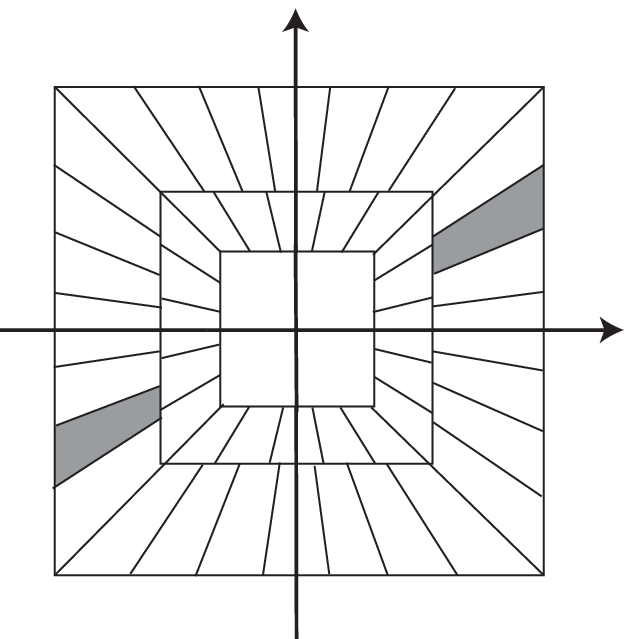}
\put(-60,-17){(b)}
\end{center}
\caption{(a) The cones $\cC_1$ -- $\cC_4$ and the centered rectangle $\cR$ in frequency domain.
(b) The tiling of the frequency domain induced by a cone-adapted shearlet system, where the (essential) support of
the Fourier transform of one shearlet generator is exemplary high-lighted.}
\label{fig:shearlet}
\end{figure}

The four cones $\cC_1$ -- $\cC_4$ are now defined by
\[
\cC_\iota = \left\{ \begin{array}{rcl}
\{(\xi_1,\xi_2) \in \bR^2 : \xi_1 \ge 1,\, |\xi_2/\xi_1| \le 1\} & : & \iota = 1,\\
\{(\xi_1,\xi_2) \in \bR^2 : \xi_2 \ge 1,\, |\xi_1/\xi_2| \le 1\} & : & \iota = 2,\\
\{(\xi_1,\xi_2) \in \bR^2 : \xi_1 \le -1,\, |\xi_2/\xi_1| \le 1\} & : & \iota = 3,\\
\{(\xi_1,\xi_2) \in \bR^2 : \xi_2 \le -1,\, |\xi_1/\xi_2| \le 1\} & : & \iota = 4,
\end{array}
\right.
\]
and a centered rectangle
\[
\cR = \{(\xi_1,\xi_2) \in \bR^2 : \|(\xi_1,\xi_2)\|_\infty < 1\}.
\]
For an illustration, we refer to Figure \ref{fig:shearlet}(a).

The rectangle $\cR$ corresponds to the low frequency content of a signal and is customarily
represented by translations of some scaling function. Anisotropy comes into play when
encoding the high frequency content of a signal which corresponds to the cones $\cC_1$ -- $\cC_4$,
where the cones $\cC_1$ and $\cC_3$ as well as $\cC_2$ and $\cC_4$ are treated separately as can be seen
in the following

\begin{definition}
\label{defi:discreteshearlets}
For some sampling constant $c > 0$, the {\em cone-adapted shearlet system} $\cSH(\phi,\psi,\tilde{\psi}\wq{;c})$ generated by
a {\em scaling function} $\phi \in L^2(\mathbb{R}^2)$ and {\em shearlets} $\psi, \tilde{\psi} \in L^2(\mathbb{R}^2)$ is
defined by
\[
\cSH(\phi,\psi,\tilde{\psi};c) = \Phi(\phi;c) \cup \Psi(\psi;c) \cup \tilde{\Psi}(\tilde{\psi};c),
\]
where
\[
\Phi(\phi;c) = \{\phi_m  = \phi(\cdot-cm) : m \in \bZ^2\},
\]
\[
\Psi(\psi;c) = \{\psi_{j,k,m} =  2^{3j/4} {\psi}({S}_{k} {A}_{2^j}\cdot-cm) :
j \ge 0, |k| \le \lceil 2^{j/2} \rceil, m \in \bZ^2 \},
\]
and
\[
\tilde{\Psi}(\tilde{\psi};c) = \{\tilde{\psi}_{j,k,m} =  2^{3j/4} \tilde{\psi}(S^T_{k} \tilde{A}_{2^j}\cdot-cm) :
j \ge 0, |k| \le \lceil 2^{j/2} \rceil, m \in \bZ^2 \}.
\]
\end{definition}
The tiling of frequency domain induced by $\cSH(\phi,\psi,\tilde{\psi};c)$ is illustrated in Figure \ref{fig:shearlet}(b).
From this illustration, the anisotropic footprints of shearlets contained in $\Psi(\psi;c)$ and $\tilde{\Psi}(\tilde{\psi};c)$
can clearly be seen. The corresponding anisotropic footprints of shearlets {\em in spatial domain} are of size $2^{-j/2}$ times
$2^{-j}$.

The reader should keep in mind that although not indicated by the notation, the functions $\phi_m$, $\psi_{j,k,m}$, and
$\tilde{\psi}_{j,k,m}$ all depend on the sampling constant $c$. For the sake of brevity, we will often write $\psi_\lambda$ and $\tilde{\psi}_{\lambda}$,
where $\lambda = (j,k,m)$ index scale, shear, and position. For later use, we further let $\Lambda_j$ and $\tilde{\Lambda}_j$ be the indexing sets
of shearlets in $\Psi(\psi;c)$ and $\tilde{\Psi}(\tilde{\psi};c)$ at scale $j$,
respectively, i.e.,
\[
\Psi(\psi;c) = \{\psi_{\lambda} : \lambda \in \Lambda_j, j=0, 1, \ldots\}
\: \mbox{ and }\:
\tilde{\Psi}(\tilde{\psi};c) = \{\tilde{\psi}_{\lambda} : \lambda \in \tilde{\Lambda}_j, j=0, 1, \ldots\}.
\]
Finally, we define
\[
\Lambda = \bigcup_{j=0}^\infty \Lambda_j \quad \mbox{and} \quad \tilde{\Lambda} = \bigcup_{j=0}^\infty \tilde{\Lambda}_j.
\]

The shearlet systems $\cSH(\phi,\psi,\tilde{\psi};c)$ have already been very well studied with respect to their
frame properties for $L^2(\RR^2)$, and we would like to refer to results in \cite{GKL06,KL07,DKST09}. It should
be mentioned that those results typically concern frame properties of $\Psi(\psi;c)$, which immediately imply
frame properties of $\tilde{\Psi}(\tilde{\psi};c)$ likewise, whereas numerous frame properties for the low-frequency part
$\Phi(\phi;c)$ can be found in the wavelet literature. Combining those leads to frame properties of $\cSH(\phi,\psi,\tilde{\psi};c)$.

Recent results in \cite{KKL10a} establish frame properties specifically for the case of spatially compactly supported shearlet
systems, i.e., shearlet systems with compactly supported generators $\phi$, $\psi$, and $\tilde{\psi}$ which lead to a
shearlet system consisting of compactly supported elements. These results give sufficient conditions for the so-called
$t_q$ conditions to be satisfied. As one class of examples with `good' frame bounds, generating shearlets $\psi$ and
$\tilde{\psi}$ were chosen to be separable, i.e., of the form $\psi_1(x_1)\cdot\psi_2(x_2)$  and $\psi_1(x_2)\cdot\psi_2(x_1)$,
respectively, where $\psi_1$ is a wavelet and $\psi_2$ a scaling function both associated with some carefully chosen (maximally
flat) low pass filter. The separability has in addition the advantage to lead to fast accompanying algorithms.

We wish to mention that there is a trade-off between {\em compactly support} of the shearlet generators, {\em tightness} of the associated frame,
and {\em separability} of the shearlet generators. The known constructions of tight shearlet frames do not use separable generators, and
these constructions can be shown to {\em not} be applicable to compactly supported generators. Tightness is difficult to obtain while
allowing for compactly supported generators, but we can gain separability, hence fast algorithmic realizations. On the other hand, when allowing
non-compactly supported generators, tightness is possible, but separability seems to be out of reach, which
makes fast algorithmic realizations very difficult.

%*******************************************************************************************
%*******************************************************************************************

\subsection{Compactly Supported Shearlet Frames for $L^2(\Omega)$}
\label{subsec:extension}

Let now $\Omega \in STAR^2(\nu,L)$ be a bounded domain as defined in Subsection \ref{subsec:ourmodel}. The main
idea to introduce a shearlet frame for $L^2(\Omega)$, preferably with compactly supported elements, is to start
with a compactly supported shearlet frame for $L^2(\RR^2)$ and apply the orthogonal projection onto $L^2(\Omega)$
to each element. To make this mathematically precise, we let $P_\Omega : L^2(\RR^2) \to L^2(\Omega)$ denote the
orthogonal projection onto $L^2(\Omega)$.

\begin{definition}
\label{defi:discreteshearlets2}
Let $\Omega \in STAR^2(\nu,L)$. For some sampling constant $c > 0$, the {\em cone-adapted shearlet system}
$\cSH_\Omega(\phi,\psi,\tilde{\psi};c)$ for $L^2(\Omega)$ generated by a {\em scaling function} $\phi \in L^2(\mathbb{R}^2)$ and
{\em shearlets} $\psi, \tilde{\psi} \in L^2(\mathbb{R}^2)$ is defined by
\[
\cSH_\Omega(\phi,\psi,\tilde{\psi};c) = P_\Omega(\Phi(\phi;c) \cup \Psi(\psi;c) \cup \tilde{\Psi}(\tilde{\psi};c)),
\]
where $\Phi(\phi;c)$, $\Psi(\psi;c)$, and $\tilde{\Psi}(\tilde{\psi};c)$ are defined as in Definition \ref{defi:discreteshearlets}.
\end{definition}

As a direct corollary from well known results in frame theory (see \cite{Chr03}), we obtain the following result,
which clarifies frame properties for systems $\cSH_\Omega(\phi,\psi,\tilde{\psi};c)$ to the extent to which they are known for
systems $\cSH(\phi,\psi,\tilde{\psi};c)$. \gggk{In the sequel, we will usually regard $\cSH_\Omega(\phi,\psi,\tilde{\psi};c)$ as
a system defined on $[0,1]^2$ -- in accordance with our model $\cE^2_{\nu,L}(\Omega)$ --
by which we simply mean extension by zero. This system will be sometimes referred to as the
{\em extension of $\cSH_\Omega(\phi,\psi,\tilde{\psi};c)$ to $[0,1]^2$.} The following result also provides frame properties
of these systems.}

\begin{proposition}
\label{prop:equalframes}
Let $c > 0$, let $\phi,\psi, \tilde{\psi} \in L^2(\mathbb{R}^2)$, and let $\Omega \in STAR^2(\nu,L)$ \wq{with positive measure}.
Then the following conditions are equivalent.
\begin{enumerate}
\item[(i)] The shearlet system $\cSH(\phi,\psi,\tilde{\psi};c)$ is a frame for $L^2(\RR^2)$ with frame bounds $A$ and $B$.
\item[(ii)] The shearlet system $\cSH_\Omega(\phi,\psi,\tilde{\psi};c)$ is a frame for $L^2(\Omega)$ with frame bounds $A$ and $B$.
\item[(iii)] \gggk{The extension of the shearlet system $\cSH_\Omega(\phi,\psi,\tilde{\psi};c)$ to $[0,1]^2$ is a frame  with frame bounds $A$ and $B$
for functions $L^2([0,1]^2)$ with compact support in $\Omega$.}
\end{enumerate}
\end{proposition}

%*******************************************************************************************
%*******************************************************************************************
%*******************************************************************************************

\section{Optimal Sparsity of Shearlets on Bounded Domains}
\label{sec:3}

We now have all ingredients to formally state the result already announced in Subsection \ref{subsec:surprise},
which shows that even with the `crude' construction of shearlets on bounded domains {\em and} the significantly
more sophisticated model for cartoon-like images on bounded domains we still obtain (almost) optimally sparse
approximations.

\begin{theorem}\label{theo:main}
Let $c > 0$, and let $\phi, \psi, \tilde{\psi} \in L^2(\R^2)$ be compactly supported. Suppose
that, in addition, for all $\xi = (\xi_1,\xi_2) \in \RR^2$, the shearlet $\psi$ satisfies
\bitem
\item[(i)] $|\hat\psi(\xi)| \le C_1 \cdot \min(1,|\xi_1|^{\alpha}) \cdot \min(1,|\xi_1|^{-\gamma}) \cdot \min(1,|\xi_2|^{-\gamma})$, and
\item[(ii)] $\left|\frac{\partial}{\partial \xi_2}\hat \psi(\xi)\right|
\le  |h(\xi_1)| \cdot \left(1+\frac{|\xi_2|}{|\xi_1|}\right)^{-\gamma}$,
\eitem
where $\alpha > 5$, $\gamma \ge 4$, $h \in L^1(\bR)$, and $C_1$ is a constant, and suppose that
the shearlet $\tilde{\psi}$ satisfies (i) and (ii) with the roles of $\xi_1$ and $\xi_2$ reversed.
Further, let $\nu > 0,L \in \Z^+$ and $\Omega \in STAR^2(\nu,L)$, and suppose that $\cSH_\Omega(\phi,\psi,\tilde{\psi};c)$
forms a frame for $L^2(\Omega)$.

Then, the \gggk{extension of the shearlet frame $\cSH_\Omega(\phi,\psi,\tilde{\psi};c)$ to $[0,1]^2$ (cf. Subsection \ref{subsec:extension})}
provides (almost) optimally sparse approximations of
functions $f \in \cE^2_{\nu,L}(\Omega)$ in the sense that there exists some $C > 0$ such that
\[
\|f-f_N\|_2^2 \leq C\cdot N^{-2} \cdot {(\log{N})}^3 \qquad \text{as} \,\, N \rightarrow \infty,
\]
where $f_N$ is the nonlinear N-term approximation obtained by choosing the N largest shearlet coefficients
of $f$.
\end{theorem}

%*******************************************************************************************
%*******************************************************************************************

\subsection{Architecture of the Proof of Theorem \ref{theo:main}}
\label{subsec:architecture}

Before delving into the proof in the following section, we present some preparation before as well
as describe the architecture of the proof for clarity purposes.

Let now $\cSH_\Omega(\phi,\psi,\tilde{\psi};c)$ satisfy the hypotheses in Theorem \ref{theo:main}, and
let $f \in \cE^2_{\nu,L}(\Omega)$. We first observe that, without loss of generality, we might
assume the scaling index $j$ to be sufficiently large, since $f$ as well as all frame elements in the
shearlet frame $\cSH_\Omega(\phi,\psi,\tilde{\psi};c)$ are compactly supported in spatial domain, hence a
finite number does not contribute to the asymptotic estimate we aim for. In particular, this means that
we do not need to take frame elements from $\Phi(\phi;c)$ into account. Also, we are allowed to restrict our
analysis to shearlets $\psi_{j,k,m}$, since the frame elements $\widetilde{\psi}_{j,k,m}$ can be
handled in a similar way.

We further observe that we can drive the analysis for the frame $\cSH(\phi,\psi,\tilde{\psi};c)$ and for
the domain $[0,1]^2$ instead, since, by hypothesis, $\Omega$ is contained in the interior of $[0,1]^2$,
we treat the boundary of $\Omega$ as a singularity curve in $[0,1]^2$, and the frame properties are
equal as shown in Proposition \ref{prop:equalframes}. In this viewpoint, the function to be sparsely
approximated vanishes on $[0,1]^2 \setminus \Omega$.

Our main concern will now be to derive appropriate estimates for the shearlet coefficients
$\{\ip{f}{\psi_{\lambda}} : \lambda \in \Lambda\}$ of $f$. Letting $|\theta(f)|_n$ denote
the $n$th largest shearlet coefficient $\ip{f}{\psi_{\lambda}}$ in absolute value and exploring
the frame property of  $\cSH(\phi,\psi,\tilde{\psi};c)$, we conclude that
\[
\|f-f_N\|_2^2 \leq \frac{1}{A}\sum_{n>N} |\theta(f)|_n^2,
\]
for any positive integer $N$, where $A$ denotes the lower frame bound of the shearlet frame $\cSH(\phi,\psi,\tilde{\psi};c)$.
Thus, for the proof of Theorem \ref{theo:main}, it suffices to show that
\begin{equation}\label{eq:upper}
\sum_{n>N} |\theta(f)|_n^2 \leq C \cdot N^{-2} \cdot  {(\log{N})}^3 \qquad \text{as} \,\, N \rightarrow \infty.
\end{equation}

%\gk{Do we say anywhere anything about whether we only consider one of the cones? \wq{Yes, on page 10.}}

To derive the anticipated estimate in \eqref{eq:upper}, for any shearlet
$\psi_{\lambda}$, we will study two separate cases:
\bitem
\item {\em Case 1}. The compact support of the shearlet $\psi_{\lambda}$ does not intersect the boundary of the set $B$ (or $\partial \Omega$), i.e.,
$\Sp(\psi_{\lambda}) \cap (\partial B \cup \partial \Omega) = \emptyset.$
\item {\em Case 2}. The compact support of the shearlet $\psi_{\lambda}$ does intersect the boundary of the set $B$ (or $\partial \Omega$), i.e.,
$\Sp(\psi_{\lambda}) \cap (\partial B \cup \partial \Omega) \ne \emptyset.$
\eitem
Notice that this exact distinction is only possible due to the spatial compact support of all shearlets
in the shearlet frame.

Case 2 will then throughout the proof be further subdivided into the situations -- which we now
do not state precisely, but just give the reader the intuition behind them:
\bitem
\item {\em Case 2a}. The support of the shearlet does intersect \wql{only one $C^2$ curve in} $\partial B \cup \partial \Omega$.
\item {\em Case 2b}. The support of the shearlet does intersect \wql{at least two $C^2$ curves in} $\partial B \cup \partial \Omega$.
\bitem
\item {\em Case 2b-1}. The support of the shearlet does intersect $\partial B \cup \partial \Omega$ in a corner point.
\item {\em Case 2b-2}. The support of the shearlet does intersect two \wql{$C^2$} curves in $\partial B \cup \partial \Omega$
simultaneously, but does not intersect a corner point.
\eitem
\eitem

%*******************************************************************************************
%*******************************************************************************************
%*******************************************************************************************

\section{Proof of Theorem \ref{theo:main}}
\label{sec:4}

In this section, we present the proof of Theorem \ref{theo:main}, following the road map outlined in
Subsection \ref{subsec:architecture}. We wish to mention that Case 1 and Case 2a are similar to the proof
of (almost) optimally sparse approximations of the class $\cE^2(\nu)$ using compactly supported shearlet frames
in \cite{KL10}. However, Case 2b differs significantly from it, since it, in particular, requires a
careful handling of the corner points of $\partial B$ and $\partial \Omega$.

In the sequel -- since we are concerned with an asymptotic estimate --
for simplicity we will often simply use $C$ as a constant although it might differ for each estimate. Also
all the results in the sequel are independent on the sampling constant $c>0$, wherefore we now fix it once
and for all.

\subsection{Case 1: The Smooth Part}

We start with Case 1, hence the smooth part. Without loss of generality, we can consider some $g \in C^2([0,1]^2)$
and estimate its shearlet coefficients. The following proposition, which is taken from \cite{KL10}, implies the
rate for optimal sparsity. Notice that the hypothesis on $\psi$ of the following result
is implied by condition (i) in Theorem \ref{theo:main}.

\begin{proposition}[\cite{KL10}]
\label{prop:main1}
Let $g \in C^2([0,1]^2)$, and let $\psi \in L^2(\RR^2)$ be compactly supported and satisfy
\[
|\hat\psi(\xi)| \le C_1 \cdot \min(1,|\xi_1|^{\alpha}) \cdot \min(1,|\xi_1|^{-\gamma}) \cdot \min(1,|\xi_2|^{-\gamma})
\:\mbox{ for all }\xi = (\xi_1,\xi_2) \in \RR^2,
\]
where $\gamma>3$, $\alpha > \gamma+2$, and $C_1$ is a constant. Then, there exists
some $C>0$ such that
\[
\sum_{n>N} |\theta(g)|_n^2 \le C \cdot N^{-2} \qquad \text{as } N \rightarrow \infty.
\]
\end{proposition}

This settles Theorem \ref{theo:main} for this case.

%*******************************************************************************************
%*******************************************************************************************

\subsection{Case 2: The Non-Smooth Part}

Next, we turn our attention to the non-smooth part, and aim to estimate the shearlet coefficients of those
shearlets whose spatial support intersects the discontinuity curve $\partial B$ or the boundary of the domain
$\Omega$. One of the main means of the proof will be the partitioning of the unit cube $[0,1]^2$ into dyadic cubes,
picking those which contain such an intersection, and estimating the associated shearlet coefficients.
For this, we first need to introduce the necessary notational concepts.

For any scale $j \ge 0$ and any grid point $p \in \ZZ^2$, we let $\cQ_{j,p}$ denote the dyadic cube defined by
\[
\cQ_{j,p} = [-2^{-j/2},2^{-j/2}]^2+2^{-j/2}p.
\]
Further, let $\cQ_j$ be the collection of those dyadic cubes $\cQ_{j,p}$ which intersect $\partial B \cup \partial \Omega$, i.e.,
\[
\cQ_j = \{\cQ_{j,p}: \cQ_{j,p} \cap (\partial B \cup \partial \Omega) \ne \emptyset, p \in \ZZ^2\}.
\]
Of interest to us is also the set of shearlet indices, which are associated with shearlets intersecting
the discontinuity curve inside some $\cQ_{j,p}  \in \cQ_{j}$, hence, for $j\ge0$ and $p \in \bZ^2$ with $\cQ_{j,p} \in \cQ_{j}$,
we will consider the index set
\[
{\Lambda_{j,p}} = \{\lambda \in \Lambda_j : \Sp(\psi_{\lambda}) \cap \cQ_{j,p} \cap (\partial B \cup \partial \Omega) \ne \emptyset\}.
\]
Finally, for $j \ge 0$, $p \in \ZZ^2$, and $0 < \eps < 1$, we define $\Lambda_{j,p}(\eps)$ to be the index set of shearlets
$\psi_\lambda$, $\lambda \in {\Lambda_{j,p}}$, such that the magnitude of the corresponding shearlet
coefficient $\langle f,\psi_\lambda \rangle$ is larger than $\eps$ and the support of $\psi_\lambda$
intersects $\cQ_{j,p}$ at the $j$th scale, i.e.,
\[
{\Lambda_{j,p}(\eps)} = \{\lambda \in {\Lambda_{j,p}} : |\langle f,\psi_{\lambda}\rangle| > \eps\},
\]
and we define {$\Lambda(\eps)$} to be the index set for shearlets so that $|\langle f,\psi_{\lambda} \rangle| > \eps$
across all scales $j$, i.e.,
\[
\Lambda(\eps) = \bigcup_{j,p} \Lambda_{j,p}(\eps).
\]
The expert reader will have noticed that in contrast to the proofs in \cite{CD04} and \cite{GL07},
which also split the domain into smaller scale boxes, we do not apply a weight function to
obtain a smooth partition of unity. In our case, this is not necessary due to the spatial compact
support of the frame elements.
Finally we set
\[
\cS_{j,p}= \bigcup_{\lambda \in \Lambda_{j,p}} \Sp(\psi_\lambda),
\]
which is contained in a cubic window of size $C\cdot 2^{-j/2}$ by $C\cdot 2^{-j/2}$, hence is of asymptotically
the same size as $\cQ_{j,p}$. As mentioned earlier, we may assume that $j$ is sufficiently large so that it is sufficient to consider the following two cases:
\begin{itemize}
\item {\em Case 2a}.
There is only
one edge curve $\Gamma_1 \subset \partial B$ (or $\partial \Omega$) which can be parameterized by $x_1 = E(x_2)$ (or $x_2 = E(x_1))$ with $E \in C^2$ inside $\cS_{j,p}$.
For any $\lambda \in \Lambda_{j,p}$, there exists some $\hat{x} = (\hat{x}_1,\hat{x}_2)
\in \cQ_{j,p} \cap \Sp(\psi_{\lambda}) \cap \Gamma_1$.
\item {\em Case 2b}. There are two edge curves $\Gamma_1,\Gamma_2 \subset \partial B$ (or $\partial \Omega$) which can be parameterized by $x_1 = E(x_2)$ (or $x_2 = E(x_1)$) with $E \in C^2$
inside $\cS_{j,p}$.
For any $\lambda \in \Lambda_{j,p}$, there exist two distinct points $\hat{x} = (\hat{x}_1,\hat{x}_2)$ and $\hat{y} = (\hat{y}_1,\hat{y}_2)$ such that
$\hat{x}\in \cQ_{j,p} \cap \Sp(\psi_{\lambda}) \cap \Gamma_1$ and $\hat{y}\in \cQ_{j,p} \cap \Sp(\psi_{\lambda}) \cap \Gamma_2$.
\end{itemize}
In the sequel, we only consider the edge curve $\partial B$ to analyze shearlet coefficients associated with the non-smooth part, since the boundary of
the domain $\Omega$ can be handled in a similar way; see also our elaboration on the fact that WLOG we can consider the approximation on $[0,1]^2$ rather than
$\Omega$ in Subsection \ref{subsec:architecture}.

%*******************************************************************************************

\subsubsection{Case 2a: The Non-Smooth Part}

This part was already studied in \cite{KL10}, where an (almost) optimally sparse approximation rate by the class of
compactly supported shearlet frames $\cSH(\phi,\psi,\tilde{\psi};c)$ under consideration was proven, and
we refer to \cite{KL10} for the precise argumentation. For intuition purposes as well as for later usage,
we though state the key estimate, which implies (almost) optimally sparse approximation for Case 2a:

\begin{proposition}[\cite{KL10}]\label{prop:main2}
Let $\psi \in L^2(\RR^2)$  be compactly supported and satisfy
the conditions (i), (ii) in Theorem \ref{theo:main} and
assume that, for any $\lambda \in \Lambda_{j,p}$,
there exists some $\hat{x}=(\hat{x}_1,\hat{x}_2) \in \cQ_{j,p} \cap \Sp(\psi_{\lambda}) \cap \partial B$.
Let $s$ be \wq{the slope\footnote{Notice that here we regard the slope of the tangent to a curve $(E(x_2),x_2)$, i.e.,
we consider $s$ of a curve $x_1 = sx_2+b$, say. For analyzing shearlets $\tilde{\psi}_{j,k,m}$, the roles of $x_1$ and $x_2$ would need to be reversed.}
of the tangent} to the edge curve $\partial B$ at $(\hat{x}_1,\hat{x}_2)$,
i.e.,
\begin{itemize}
\item $s = E'(\hat{x}_2)$, if $\partial B$ is parameterized by $x_1 = E(x_2)$ with $E \in C^2$ in $\cS_{j,p}$,
\item $s = ({E'(\hat{x}_1)})^{-1}$, if $\partial B$ is parameterized by $x_2 = E(x_1)$ with $E \in C^2$ in $\cS_{j,p}$, and
\item $s = \infty$, if $\partial B$ is parameterized by $x_2 = E(x_1)$ with $E'(\hat{x}_1) = 0$ and $E \in C^2$ in $\cS_{j,p}$.
\end{itemize}
Then, there exists some $C > 0$ such that
\begin{equation}\label{eq:estimate1}
|\langle f,\psi_{\lambda}\rangle| \leq C \cdot 2^{-\frac{9}{4}j}, \qquad \text{if }|s| > \wq{\frac{3}{2}} \text{ or } |s| = \infty,
\end{equation}
and
\begin{equation}\label{eq:estimate2}
|\langle f,\psi_{\lambda}\rangle| \leq C \cdot \frac{2^{-\frac{3}{4}j}}{|\wq{k+}2^{j/2}s|^3},  \qquad \text{if }|s| \leq \wq{3}.
\end{equation}
\end{proposition}

Similar estimates with
$\partial B$ substituted by $\partial \Omega$ hold if, for any $\lambda \in \Lambda_{j,p}$,
there exists some $\hat{x}=(\hat{x}_1,\hat{x}_2) \in \cQ_{j,p} \cap \Sp(\psi_{\lambda}) \cap \partial \Omega$.

%*******************************************************************************************

\subsubsection{Case 2b: The Non-Smooth Part}

Letting $\eps > 0$, our goal will now be to first estimate $|\Lambda_{j,p}(\eps)|$ and, based on this, derive
an estimate for $|\Lambda(\eps)|$. WLOG we might assume $\|\psi\|_1 \le 1$, which implies
\[
|\langle f,\psi_{\lambda}\rangle| \le 2^{-\frac{3}{4}j}.
\]
Hence, for estimating $|\Lambda_{j,p}(\eps)|$, it is sufficient to restrict our attention to scales $j \le \frac{4}{3}\log_2(\eps^{-1})$.

\begin{figure}[ht]
\begin{center}
\includegraphics[height=1.4in]{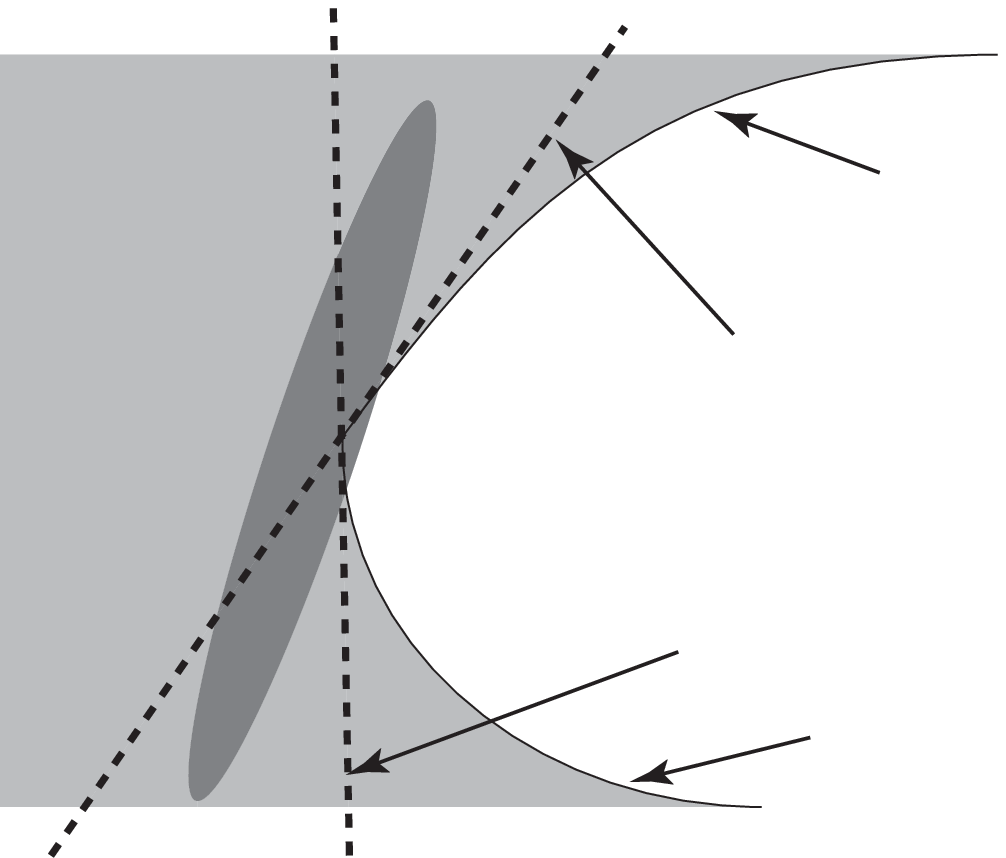}
\put(-60,44){\footnotesize{$B_1$}}
\put(-110,70){\footnotesize{$B_0$}}
\put(-33,28){\footnotesize{$T_2$}}
\put(-30,55){\footnotesize{$T_1$}}
\put(-19,13){\footnotesize{$\Gamma_2$}}
\put(-12,74){\footnotesize{$\Gamma_1$}}
\put(-70,-17){(a)}
\hspace*{1cm}
\includegraphics[height=1.4in]{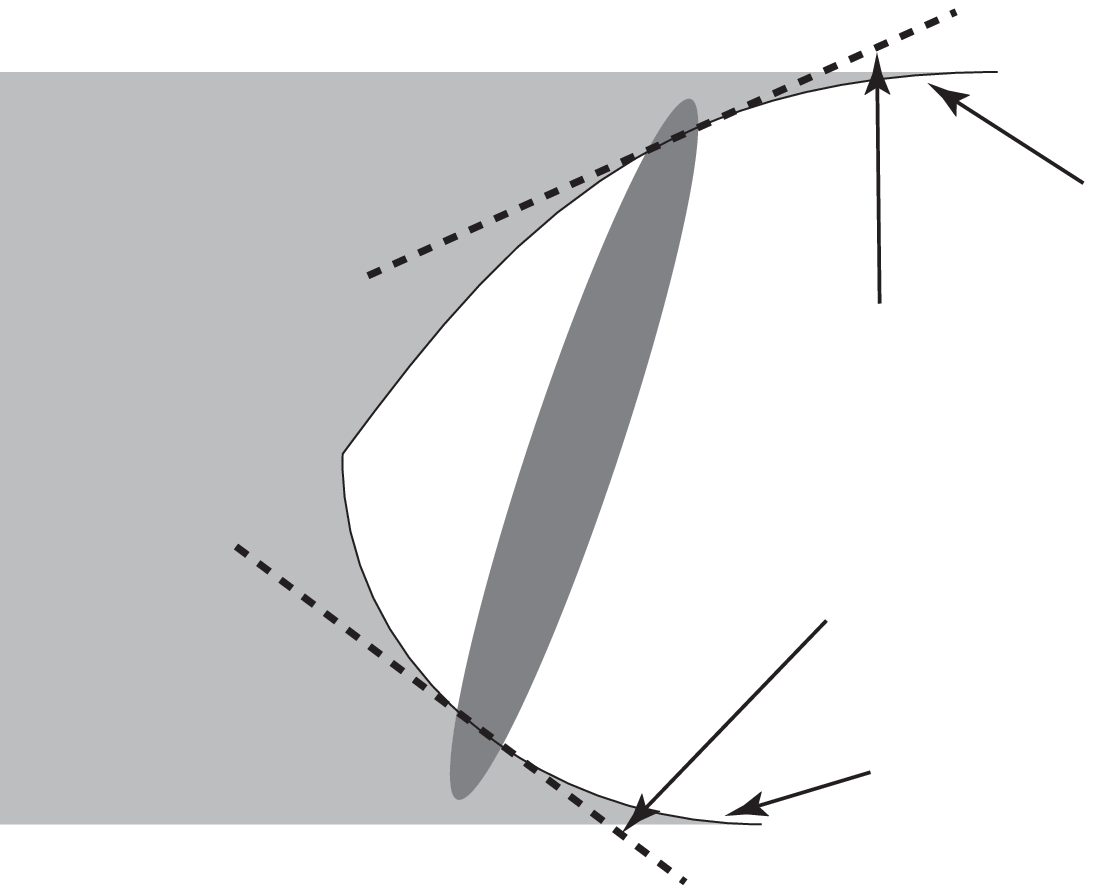}
\put(-80,44){\footnotesize{$B_1$}}
\put(-110,70){\footnotesize{$B_0$}}
\put(-26,32){\footnotesize{$T_2$}}
\put(-27,58){\footnotesize{$T_1$}}
\put(-19,13){\footnotesize{$\Gamma_2$}}
\put(1,72){\footnotesize{$\Gamma_1$}}
\put(-70,-17){(b)}
\end{center}
\caption{(a) A shearlet $\psi_{\lambda}$ intersecting a corner point where two edge curves $\Gamma_1$ and $\Gamma_2$ meet.
$T_1$ and $T_2$ are tangents to the edge curves $\Gamma_1$ and $\Gamma_2$ in this corner point.
(b) A shearlet $\psi_{\lambda}$ intersecting two edge curves $\Gamma_1$ and $\Gamma_2$ which are a part of the
boundary of sets $B_0$ and $B_1$. $T_1$ and $T_2$ are tangents to the edge curves $\Gamma_1$ and $\Gamma_2$ in points
contained in the support of $\psi_{\lambda}$.
}
\label{fig:shear}
\end{figure}

As already announced before, we now split Case 2b into the following two subcases:
\begin{itemize}
\item {\em Case 2b-1}. The shearlet $\psi_{\lambda}$ intersects a corner point, in which two $C^2$ curves $\Gamma_1$ and $\Gamma_2$, say, meet
(see Figure \ref{fig:shear} (a)).
\item {\em Case 2b-2}. The shearlet $\psi_{\lambda}$ intersects two edge curves $\Gamma_1$ and $\Gamma_2$, say, simultaneously, but it does not intersect
a corner point (see Figure \ref{fig:shear} (b)).
\end{itemize}

{\em Case 2b-1}. We first consider {\em Case 2b-1}. In this case, by a counting argument, it follows that
\[
|\Lambda_{j,p}(\eps)| \leq C \cdot 2^{j/2}.
\]
Since there are only finitely many corner points with its number not depending on scale $j \ge 0$, we have
\[
|\Lambda(\eps)| \leq C \cdot \sum_{j=0}^{\frac{4}{3}\log_2{(\eps^{-1})}} 2^{j/2} \leq C \cdot \eps^{-\frac{2}{3}}.
\]
The value $\eps>0$ can be written as a function of the total number $N$ of coefficients, which yields
$\eps(N) \leq C \cdot N^{-\frac{3}{2}}$. This implies that
\[
\sum_{n>N} |\theta(f)|^2_n \leq C \cdot N^{-2},
\]
and the optimal sparse approximation rate is proven for {\em Case 2b-1}.

{\em Case 2b-2}. Next, we consider {\em Case 2b-2}. In this case, WLOG, we might assume that,
for any $\lambda \in \Lambda_{j,p}$, there exist two distinct points $\hat{x} = (\hat{x}_1,\hat{x}_2),
\hat{y} = (\hat{y}_1,\hat{y}_2)$ such that $\hat{x}\in \cQ_{j,p} \cap \Sp(\psi_{\lambda}) \cap \Gamma_1$ and
$\hat{y}\in \cQ_{j,p} \cap \Sp(\psi_{\lambda}) \cap \Gamma_2$, and the two edge curves $\Gamma_1$ and $\Gamma_2$
are parameterized by $x_1 = E(x_2)$ (or $x_2 = E(x_1)$) with $E \in C^2$ inside $S_{j,p}$.
We can then write the function $f \in \cE^2_{\nu,L}(\Omega)$ as
\[
f_0\chi_{B_0}+f_1\chi_{B_1} = (f_0-f_1)\chi_{B_0}+f_1 \qquad \text{on } \cS_{j,p},
\]
where $f_0,f_1 \in C^2([0,1]^2)$ and $B_0,B_1$ are two disjoint subsets of $[0,1]^2$ (see Figure \ref{fig:shear}).
By Proposition \ref{prop:main1}, the rate for optimal sparse approximation is achieved for the smooth part
$f_1$. Thus, it is sufficient to consider $f = g_0\chi_{B_0}$ with $g_0 = f_0-f_1 \in C^2([0,1]^2)$.

Assume now that the tangents to the edge curves $\Gamma_1$ and $\Gamma_2$ at the points $\hat{x}$ and $\hat{y}$ are
given by the equations
\[
T_1 \: : \: x_1+a_1 = s_1(x_2+b_1) \quad \text{and} \quad T_2 \: : \: x_1+a_2 = s_2(x_2+b_2),
\]
respectively, i.e., $s_1$ and $s_2$ are \wq{the slopes of the tangents} to the edge curves $\Gamma_1$
and $\Gamma_2$ at $\hat{x}$ and $\hat{y}$, respectively. If the curve $\Gamma_i$, $i=1,2$, is parameterized by $x_2 = E(x_1)$
with $E'(\hat{x}_1) = 0$, we let $s_i = \infty$ so that the tangent is given by $x_2 = -b_i$ in this case.

Next, for fixed scale $j$ and shear index $k$,
let $N^1_{j,k}(Q_{j,p})$ denote the number of shearlets $\psi_{\lambda}$ intersecting $\Gamma_1$ in $\cQ_{j,p}$, i.e.,
\[
N^1_{j,k}(Q_{j,p}) = |\{\lambda = (j,k,m) : \cQ_{j,p} \cap \Sp(\psi_{\lambda}) \cap \Gamma_1 \neq \emptyset\}|,
\]
let $N^2_{j,k}(Q_{j,p})$ denote the number of shearlets $\psi_{\lambda}$ intersecting $\Gamma_2$ in $\cQ_{j,p}$,
i.e.,
\[
N^2_{j,k}(Q_{j,p}) = |\{\lambda = (j,k,m) :\cQ_{j,p} \cap \Sp(\psi_{\lambda}) \cap \Gamma_2 \neq \emptyset \}|,
\]
and let $N_{j,k}(Q_{j,p})$ denote the number of shearlets $\psi_{\lambda}$ intersecting $\Gamma_1$ and $\Gamma_2$ in $\cQ_{j,p}$,
i.e.,
\[
N_{j,k}(Q_{j,p}) = |\{\lambda = (j,k,m) :\cQ_{j,p} \cap \Sp(\psi_{\lambda}) \cap \Gamma_1 \neq \emptyset\mbox{ and }\cQ_{j,p} \cap \Sp(\psi_{\lambda}) \cap \Gamma_2 \neq \emptyset \}|.
\]
\wq{Then}
\beq \label{eq:N1}
N_{j,k}(Q_{j,p}) \leq \min{(N^1_{j,k}(Q_{j,p}),N^2_{j,k}(Q_{j,p}))}.
\eeq
By a counting argument, there exists some $C > 0$ such that
\beq \label{eq:N2}
N^i_{j,k}(Q_{j,p}) \leq C \cdot 2^{j/2} \qquad \mbox{for } i=1,2,
\eeq
and the form of $\Sp(\psi_{\lambda})$ implies
\beq \label{eq:N3}
N^i_{j,k}(Q_{j,p}) \leq C \cdot (|2^{j/2}s_i\wq{+k}|+1)\qquad \mbox{for } i=1,2.
\eeq

We now subdivide into three subcases, namely, $|s_1|,|s_2| \leq 2$, and $|s_1| \leq 2, |s_2| > 2$ (or vice versa), and $|s_1|,|s_2|>2$,
and show in each case the (almost) optimal sparse approximation rate claimed in Theorem \ref{theo:main}. This then finishes
the proof.

{\em Subcase $|s_1|,|s_2| \leq 2$.} In this case, \eqref{eq:N1} and \eqref{eq:N3} yield
\[
N_{j,k}(Q_{j,p}) \leq C \cdot \min{(|2^{j/2}s_1\wq{+k}|+1,|2^{j/2}s_2\wq{+k}|+1)}.
\]
We first show independence on the values of $s_1$ and $s_2$ within the interval $[-2,2]$. For this, let $s$ and $s'$ be the
\wq{slopes of the tangents to the edge curve $\Gamma_1$ (or $\Gamma_2$) at $t \in Q_{j,p} \cap \Sp({\psi_{\lambda}})$ and $t' \in Q_{j,p} \cap \Sp({\psi_{\lambda'}})$,
respectively, with $s \in [-2,2]$. Since $\Gamma_1$ (or $\Gamma_2$) is $C^2$, we have $|s-s'| \leq C \cdot 2^{-j/2}$, and hence
\[
|2^{j/2}s'\wq{+k}| \leq C \cdot (|2^{j/2}s\wq{+k}|+1).
\]
This implies that the estimate for $N_{j,k}(Q_{j,p})$ asymptotically remains the same, independent of the values of $s_1$ and $s_2$.
Further, we may assume $s' \in [-3,3]$ for $s \in [-2,2]$, since a scaling index $j$ can be chosen such that $|s-s'|$ is sufficiently small.
Therefore, one can apply inequality \eqref{eq:estimate2} from Proposition \ref{prop:main2} for both points $t$ and $t'$.
In fact, it can be easily checked that one can use \eqref{eq:estimate2} with the slope $s$ instead of $s'$ for the point $t'$ (or vice versa);
this replacement will not change the asymptotic estimates which we will derive.
Thus, we might from now on use universal values for the slopes $s_1$ and $s_2$ at each
point in $Q_{j,p}$.}

Now using \eqref{eq:estimate2} from Proposition \ref{prop:main2}, we have
\begin{equation}\label{eq:estimate3}
|\langle f,\psi_{j,k,m}\rangle| \leq C \cdot \max \Bigl( \frac{2^{-\frac{3}{4}j}}{|2^{j/2}s_1\wq{+k}|^3},\frac{2^{-\frac{3}{4}j}}{|2^{j/2}s_2\wq{+k}|^3}\Bigr).
\end{equation}
Since $\frac{2^{-\frac{3}{4}j}}{|2^{j/2}s_i\wq{+k}|^3}>\eps$ implies
\[
|2^{j/2}s_i\wq{+k}| < \eps^{-\frac{1}{3}}2^{-\frac{1}{4}j} \qquad \text{for } i=1,2,
\]
the estimate \eqref{eq:estimate3} yields
\begin{eqnarray}\label{eq:estimate4}
|\Lambda_{j,p}(\eps)| &\leq& C\cdot \sum_{k \in K^1_j(\eps) \cup K^2_j(\eps)} \min (|2^{j/2}s_1\wq{+k}|+1,|2^{j/2}s_2\wq{+k}|+1) \nonumber \\
&\leq& C\cdot \sum_{i=1}^{2}\sum_{k \in K^i_j(\eps)}(|2^{j/2}s_{i}\wq{+k}|+1) \nonumber \\
&\leq& C\cdot \wq{(\eps^{-\frac{1}{3}}2^{-\frac{1}{4}j}+1)^2},
\end{eqnarray}
where
\[
K^i_j(\eps) = \{k \in \Z : |2^{j/2}s_i\wq{+k}|<\wq{C} \cdot \eps^{-\frac{1}{3}}2^{-\frac{1}{4}j}\} \quad \text{for} \,\,\,\, i = 1,2.
\]
By the hypothesis for Case 2b-2, we have $|Q_j| \le C$, where the constant $C$ is independent of scale $j \ge 0$.
Therefore, continuing \eqref{eq:estimate4},
\[
|\Lambda(\eps)| \leq C \cdot\sum_{j=0}^{\frac{4}{3}\log_2 (\eps^{-1})} |\Lambda_{j,p}(\eps)|
\leq C\cdot \eps^{-\frac{2}{3}}.
\]
This allows us to write $\eps>0$ as a function of the total number of coefficients $N$, which gives
\[
\eps(N) \leq C \cdot N^{-\frac{3}{2}}.
\]
Thus
\begin{equation}\label{eq:estimate5}
\sum_{n>N}|\theta(f)|_n^2 \leq C \cdot N^{-2},
\end{equation}
which is the rate we seeked.

{\em Subcase $|s_1| \leq 2$ and $|s_2| > 2$ or vice versa.} In this case, \eqref{eq:N1}--\eqref{eq:N3} yield
\[
N_{j,k}(Q_{j,p}) \leq C \cdot \min{(|2^{j/2}s_1\wq{+k}|+1,2^{j/2})}.
\]
Again utilizing the fact that the edge curves are $C^2$, and using similar arguments as in the first subcase, WLOG we
can conclude that the slopes $s_1,s_2$ at each point in $Q_{j,p}$ are greater than $\frac 32$.

Now, exploiting inequalities \eqref{eq:estimate1} and \eqref{eq:estimate2} from Proposition \ref{prop:main2}, we have
\begin{equation}\label{eq:estimate6}
|\langle f,\psi_{j,k,m}\rangle| \leq C \cdot \max \Bigl( \frac{2^{-\frac{3}{4}j}}{|2^{j/2}s_1\wq{+k}|^3},2^{-\frac{9}{4}j} \Bigr).
\end{equation}
Since $\frac{2^{-\frac{3}{4}j}}{|2^{j/2}s_1\wq{+k}|^3} > \eps$ implies
\[
|2^{j/2}s_1\wq{+k}| < \eps^{-\frac{1}{3}}2^{-\frac{1}{4}j},
\]
and $2^{-\frac{9}{4}j} > \eps$ implies
\[
j \leq \frac{4}{9}\log_2 (\eps^{-1}),
\]
it follows from \eqref{eq:estimate6}, that
\begin{eqnarray*}
|\Lambda(\eps)|
&\leq& C\cdot \Bigl(\sum_{j=0}^{\frac{4}{3}\log_2(\eps^{-1})}\sum_{k \in K^1_j(\eps)} ( |2^{j/2}s_1\wq{+k}|+1 ) + \sum_{j=0}^{\wq{\frac{4}{9}}\log_2(\eps^{-1})}2^{j/2}\Bigr) \\
&\leq& C\cdot \Bigl(\sum_{j=0}^{\frac{4}{3}\log_2(\eps^{-1})}(\eps^{-\frac{1}{3}}2^{-j/4}+1)^2 + \sum_{j=0}^{\frac{4}{9}\log_2(\eps^{-1})}2^{j/2}\Bigr) \\
&\leq& C\cdot \eps^{-\frac{2}{3}}.
\end{eqnarray*}
The value $\eps>0$ can now be written as a function of the total number of coefficients $N$, which gives
\[
\eps(N) \leq C \cdot N^{-\frac{3}{2}}.
\]
Thus, we derive again the seeked rate
\[
\sum_{n>N}|\theta(f)|_n^2 \leq C \cdot N^{-2}.
\]

{\em Subcase $|s_1|>2$ and  $|s_2|>2$.} In this case, \eqref{eq:N1} and \eqref{eq:N2} yield
\[
N_{j,k}(Q_{j,p}) \leq C \cdot 2^{j/2}.
\]
Following similar arguments as before, we again derive the seeked rate \eqref{eq:estimate5}.

%*******************************************************************************************
%*******************************************************************************************
%*******************************************************************************************

\section{Discussion}
\label{sec:5}

A variety of applications are concerned with efficient encoding of 2D functions defined on non-rectangular domains
exhibiting curvilinear discontinuities, such as, e.g., a typical solution of a transport dominated partial differential
equation. As an answer to this problem, our main result, Theorem \ref{theo:main}, shows that compactly supported
shearlets satisfying some weak decay and smoothness conditions, when orthogonally projected onto a given domain bounded
by a piecewise $C^2$ curve, provide (almost) optimally sparse approximations of functions which are $C^2$ apart from a
piecewise $C^2$ discontinuity curve. In this model the boundary curve is treated as a discontinuity curve.

Analyzing the proof of Theorem \ref{theo:main}, it becomes evident that the presented optimal sparse
approximation result for functions in $\cE^2_{\nu,L}(\Omega)$ generalizes to an even more encompassing
model, which does contain multiple piecewise $C^2$ possibly intersecting discontinuity curves separating
$C^2$ regions in the bounded domain $\Omega$.

\smallskip

In some applications it is though of importance to avoid discontinuities at the boundary of the domain. Tackling
this question requires further studies to carefully design shearlets near the boundary, and this will be one
of our objective for the future.

%*******************************************************************************************
%*******************************************************************************************
%*******************************************************************************************

\begin{acknowledgement}
Both authors would like to thank Wolfgang Dahmen, Philipp Grohs, Chunyan Huang, Demetrio Labate, Christoph Schwab,
and Gerrit Welper for various discussions on related topics, and Jakob Lemvig for careful reading of an earlier
version of this paper. They acknowledge support from DFG Grant SPP-1324,
KU 1446/13. The first author also acknowledges support from DFG Grant KU 1446/14.
\end{acknowledgement}

\end{document}